\magnification=\magstep1
\nopagenumbers
\headline={\ifnum\pageno=1\hfil\else\hss--\ \tenrm\folio\ --\hss\fi}
\def\a{\alpha}
\def\aa{a^\delta_\omega}
\def\cl{\mathop{\rm cl}}
\def\d{\delta}

\def\e{\emptyset}
\def\forces{\parallel\!\!\!-}

\def\k{\kappa}
\def\l{\lambda}
\def\le{\langle}
\def\L{L^\infty(Y,{\cal B},\nu)}
\def\m{\medskip}
\def\n{\noindent}
\def\o{\omega}
\def\ovr{\overline}
\def\r{\rangle}
\def\R{{\bf R}}
\def\s{\sigma}
\def\shelah{{\bf[S 83]}}
\def\sq{\subseteq}

\def\v{\vskip .1cm}
\def\endproof{\ \ \vrule height6pt width6pt depth1pt}
\hfill version of 92/01/03
\m
\centerline{\bf Linear liftings for non-complete probability spaces}
\vskip 1cm
\centerline{Maxim R. Burke\footnote{$^{1)}$}{Research supported by
UPEI Senate Grant no 602101, by the Research Institute for Mathematical
Sciences at Bar-Ilan University, and by the Landau Center for Mathematical
Research in
Analysis (supported by the Minerva Foundation).  The author would like
to thank the organizers of the Winter Institute on the Set Theory of the Reals
for their hospitality while part of this research was being carried out.}}
\vskip .2cm
\centerline{and}
\vskip .2cm
\centerline{Saharon Shelah\footnote{$^{2)}$}{Partially supported by the
Foundation for Basic Research of the Israel Academy of Science. Publication
number 437.}}
\vskip 3cm
{\bf Abstract}:  We show that it is consistent with ZFC that $\L$ has no linear
lifting for many non-complete probability spaces $(Y,{\cal B},\nu)$, in
particular for $Y=[0,1]^A$, ${\cal B}=$ Borel subsets of $Y$, $\nu=$ usual
Radon measure on ${\cal B}$.
\vskip 5cm
AMS Subject Classification 1980 (1985 revision):

Primary: 28A51 \qquad Secondary: 03E35
\vfill\eject
\noindent{\bf\S 1  Introduction}
\vskip .2cm
In \shelah\ the second author showed that it is consistent that Lebesgue measure
on $[0,1]$ has no Borel lifting.  This argument was generalized in {\bf[J 89]}
and {\bf[BJ 89]} to produce a model where there is no lifting $\rho$ for the
usual product measure on $[0,1]^A$ such that for each measurable set $E$,
$\rho(E)=E'\times[0,1]^{A-B}$ where $B\sq A$ is countable and $E'\sq[0,1]^B$ is
projective. In particular $[0,1]^A$ has no Baire lifting. 
The approach taken there did not shed any light on the question of 
whether one can produce in ZFC a
Borel lifting for $[0,1]^A$ when $A$ is uncountable. In this paper we show that
this is not possible.  D.H. Fremlin suggested the use of linear liftings for
this purpose. The technique is a modification of the one used in
\shelah. We assume that the reader is familiar with \shelah. Most
definitions which we need are given below.  See {\bf[IT 69]} for more details
concerning liftings.
\m
\n{\bf1.1 Definitions} 
\m
\item{1.} If $(Y,{\cal B},\nu)$ is any probability space (not necessarily 
complete) 
then as usual we say that $f\colon Y\to\R$ is {\it measurable} if
$f^{-1}(a,b)\in{\cal B}$ for every rational interval $(a,b)\sq\R$. 
\v
\item{2.} $\L=\{f\in\R^Y: f$ is bounded and measurable $\}$.
\v
\item{3.} $\rho\colon\L\to\L$ is a {\it linear lifting} if for all $f,g\in\L$ 
and all  $x,y\in\R$,
\v
\itemitem{(a)} $f=g$ a.e. implies $\rho(f)=\rho(g)$ (everywhere).
\v
\itemitem{(b)} $\rho(f)=f$ a.e.
\v
\itemitem{(c)} $\rho(xf+yg)=x\rho(f)+y\rho(g)$.
\v
\itemitem{(d)} $\rho(1)=1$ where $1$ is the constant function with value 1.
\v
\itemitem{(e)} $f\geq0$ a.e. implies $\rho(f)\geq0$.
\v
\item{4.} $\rho\colon\L\to\L$ is a {\it lifting} if $\rho$ is a linear lifting
and
$\rho(fg)=\rho(f)\rho(g)$ for all $f,g\in\L$. In this case $\rho$ corresponds
in a canonical way 
to a lifting for the measure algebra of $(Y,{\cal B},\nu)$. See {\bf[IT
69]}.
\v
\item{5.} When $\rho$ is a linear lifting for $\L$ and $E\in{\cal B}$, we will 
write $\rho(E)$ instead of $\rho(\chi_E)$, where $\chi_E$ is the characteristic
function of the set $E$.
\v
\item{6.} For sequences of real numbers, we will use the expressions {\it 
increasing}
and {\it decreasing} to mean {\it strictly} increasing and {\it strictly}
decreasing, respectively.
\v
\item{7.} For real numbers $a\not=b$, $(a,b)$ will denote $\{x\in\R:a<x<b\}$ if
$a<b$, and $\{x\in\R:b<x<a\}$ if $b<a$.
\m
We will prove the following theorem:
\m
\n{\bf1.2 Theorem} {\it The following is consistent with ZFC\/:
\m
Let $\Sigma=$ the $\s$-algebra of Borel subsets of $[0,1]$, $\mu=$
Lebesgue measure on $\Sigma$.

Then $L^\infty([0,1],\Sigma,\mu)$ has no linear lifting.}
\vfill\eject
\n{\bf1.3 Corollary} {\it The following is consistent with ZFC\/:
\m
Suppose that
\v
{\rm1.} $(Y,{\cal B},\nu)$ is a probability space (not necessarily
complete), 
\v
{\rm2.} There is a measurable inverse-measure-preserving function
$\varphi\colon Y\to [0,1]$,
\v
{\rm3.} There is a {\rm Borel disintegration} of $\nu$, i.e., there is a family
$\le\nu_x:x\in[0,1]\r$ of probability measures on ${\cal B}$ such that for each
$g\in\L$, the function $x\mapsto\int g\,d\nu_x$ is Borel measurable and equal
a.e. to $E(g|\varphi^{-1}(\Sigma))$. (Here $E(\cdot)$ is the conditional
expectation operator.)
\m
\n Then $L^\infty(Y,{\cal B},\nu)$ has no linear lifting. 
(In particular $(Y,{\cal B},\nu)$ has no lifting.)}
\m
\n[Proof: If $\rho$ is a linear lifting for $\L$, then $\bar\rho$ is
a linear lifting for $L^\infty([0,1],\Sigma,\mu)$
where $\bar\rho(f)(x)=\int\rho(f\circ\varphi)\,d\nu_x$. [For a.a. $x$ we
have
$\bar\rho(f)(x)=E(f\circ\varphi|\varphi^{-1}(\Sigma))(x)=f(x)$.]]
\m
\n{\bf1.4 Examples} 
\m
\item{1.} $Y=[0,1]^A$, ${\cal B}=$ Borel subsets of $[0,1]^A$, $\nu=$ usual 
Radon product measure on ${\cal B}$.
\m
\item{2.} $Y=\{0,1\}^A$, ${\cal B}=$ Borel subsets of $\{0,1\}^A$, $\nu=$ usual 
Haar measure on ${\cal B}$.
\m
\item{3.} $(Z,{\cal C},\l)$ is any probability space, $Y=[0,1]\times Z$, 
${\cal B}=$ the $\s$-algebra
generated by the rectangles $E\times F$, $E\in\Sigma$, $F\in{\cal C}$, 
$\nu=$ the usual product measure on ${\cal B}$.
\m
Note that the third hypothesis of the corollary is needed. To see this,
consider the hyperstonian
space $(Y,{\cal B},\nu)$ of [0,1] and the canonical projection $\varphi\colon
Y\to[0,1]$. We know that $(Y,{\cal B},\nu)$ has a lifting (even a continuous
lifting). (See {\bf[F 89]}.) However in the model which we will construct, 
none of the disintegrations of $\nu$ will be Borel, so there is no
contradiction.
\m
\n{\bf1.5 Problem} Is it consistent with ZFC that there is a
translation-invariant linear lifting for $L^\infty([0,1),\Sigma,\mu)$?
($\rho$ is {\it translation invariant} if $\rho(f_a)(x)=\rho(f)(a+x)$, where
$f_a(y)=f(a+y)$ (all additions are mod 1), for $a,x,y\in[0,1)$,
$f\in L^\infty([0,1),\Sigma,\mu)$.)
\m
\n {\bf\S2 Proof of theorem 1.2} 
\m
Let $L^\infty$ stand for
$L^\infty([0,1],\Sigma,\mu)$.
\m
Assume $V=L$. As in \shelah\ (the technique is explained in {\bf[S 82]}), we use
an oracle-cc iteration of length 
$\aleph_2$, and it will suffice to prove the following lemma.
\m
\n{\bf2.1 Main Lemma} {\it Let $\ovr{M}$ be an $\aleph_1$-oracle and let $\rho$
be a linear
lifting of $L^\infty$. Then there is a forcing notion $P$ satisfying the
$\ovr{M}$-cc and a $P$-name $\dot X$ of an open set such that for every
$G\sq P\times Q$ generic over $V$ (where $Q$ is Cohen forcing), there is no
Borel function $h$ in $V[G]$ such that 

{\rm(a)} $h=\chi_{\dot X[G]}$ a.e.

{\rm(b)} for every $g\in(L^\infty)^V$, if $g\leq\chi_{\dot X[G]}$ a.e. then 
$\rho(g)\leq h$.

{\rm(c)} for every $g\in(L^\infty)^V$, if $\chi_{\dot X[G]}\leq g$ a.e. then 
$h\leq \rho(g)$.}
\m
\n{\bf2.2 Proof of the main lemma} 
\m
Let ${\cal S}$ denote the set of triples 
$$\bar a=\bigl(\le a_{0i}:i<\o\r,\le a_{1i}:i<\o\r,a_\o\bigr)$$ 
such that the
$a_{ji}$ are rational numbers in $(0,1)$ ($j<2$, $i<\o$), $a_\o$ is irrational,
$\le a_{0i}:i<\o\r$ is an increasing sequence converging to $a_\o$ and $\le
a_{1i}:i<\o\r$ is a decreasing sequence converging to $a_\o$. 
\v
Define a partial order $P=P(\le\bar a^\a:\a<\beta\r)$ where $\beta\leq\o_1$,
$\bar a^\alpha\in{\cal S}$, and the numbers $a^\a_\o$ are pairwise distinct, as
follows: $p\in P$ {\it iff\/} the following conditions hold:
\v
\item{(a)} $p=(U_p,f_p)$, where $U_p$ is an open subset of $(0,1)$, $\cl(U_p)$ 
has measure $<1/2$, and $f_p\colon U_p\to\{0,1\}$.
\v
\item{(b)} There is a finite sequence of rational numbers 
$0=b_0<b_1<\dots<b_n=1$ such
that $U_p=\bigcup_{\ell=0}^{n-1}I_\ell$, $\cl(I_\ell)\sq(b_\ell,b_{\ell+1})$.
\v
\item{(c)} $I_\ell$ is either a rational interval, in which case $f_p|I_\ell$ 
is constant, or there are $\a<\beta$ and $n(\ell)<\o$ such that 
$$I_\ell=\bigcup_{j<2}\ \bigcup_{n(\ell)\leq m<\o}(a^\a_{j,2m},a^\a_{j,2m+1})$$ 
and
$f_p|(a^\a_{j,4m+2k},a^\a_{j,4m+2k+1})$ is identically equal to $k$, ($j<2$,
$n(\ell)\leq2m+k$, $m<\o$, $k<2$).
\m
The order on $P$ is: $p\leq q$ if and only if $U_p\sq U_q$, $f_p\sq f_q$, and
$\cl(U_p)\cap U_q=U_p$.

Let $\dot X$ be a $P$-name for $\bigcup\{(a,b): (a,b)$ is a rational interval
$\sq(0,1)$ and for some $p\in G_P$, $(a,b)\sq U_p$ and $f_p|(a,b)$ is
identically zero$\}$. 
\m
As in \shelah, the main lemma will follow if we prove the following claim.
\m
\n{\bf2.3 Main Claim} {\it Let $P_\d=P(\le\bar a^\a:\a<\d\r)$, $\d<\o_1$ be 
given, as
well as a countable $M_\d$, $P_\d\in M_\d$, a condition $(p^*,r^*)\in
P_\d\times Q$ and a $P_\d\times Q$-name $\tau$ for a code for a member of
$L^\infty$. (We shall identify Borel functions and their codes. This should not
cause any confusion.)
Then we can find $\bar a^\d\in{\cal S}$ such that, letting
$P_{\d+1}=P(\le\bar a^\a:\a\leq\d\r)$, the following conditions hold:
\m
{\rm(A)} Every predense subset of $P_\d$ which belongs to $M_\d$ is a predense
subset of $P_{\d+1}$.
\v
{\rm(B)} There is a condition $(p',r')\in P_{\d+1}\times Q$ such that
$(p^*,r^*)\leq(p',r')$ and one of the following two conditions holds for some 
$n$:
\v
\settabs\+\indent\quad&{\rm(B1)} $(p',r')\forces_{P_{\d+1}\times Q}$ {\rm``}&\cr
\+&{\rm(B1)} $(p',r')\forces_{P_{\d+1}\times Q}$ {\rm``}&\ $\tau(\aa)\geq1/2$\cr
\+&&and
$\rho(\bigcup_{j<2}\bigcup_{n\leq m<\o}(a^\d_{j,4m+2},a^\d_{j,4m+3}))(\aa)
\geq3/4$\cr 
\+&&and $\dot X\cap
\bigcup_{j<2}\bigcup_{n\leq m<\o}(a^\d_{j,4m+2},a^\d_{j,4m+3})=\e$.{\rm''}\cr
\v
\n or 
\+&{\rm(B2)} $(p',r')\forces_{P_{\d+1}\times Q}$ {\rm``}&\ $\tau(\aa)\leq1/2$\cr
\+&&and
$\rho(\bigcup_{j<2}\bigcup_{n\leq m<\o}(a^\d_{j,4m},a^\d_{j,4m+1}))(\aa)
\geq3/4$\cr 
\+&&and 
$\bigcup_{j<2}\bigcup_{n\leq m<\o}(a^\d_{j,4m},a^\d_{j,4m+1})\sq\dot X$.{\rm''}
\cr}
\m
\n[The proof of the Main Lemma is a bookkeeping argument using the Main Claim.
$P$ is obtained, in the notation of the Main Claim, as
$P=\bigcup_{\d<\o_1}P_\d$, and the bookkeeping is needed to ensure that all
triples $(p^*,r^*,\tau)$ are considered in the construction, where
$(p^*,r^*)\in P\times Q$ and $\tau$ is a $P\times Q$-name for a code of a Borel
function.  If there were an $h$ contradicting the Main Lemma, then there would
be a $P\times Q$-name $\tau$ for $h$ and a condition $(p^*,r^*)\in P\times Q$
forcing that $\tau$ satisfies (a), (b), (c) of the Main Lemma. But then
condition (B) of the Main Claim gives a contradiction. For more details of such
oracle-cc arguments see pp. 114ff of {\bf[S 82]}.]
\m
\n{\bf2.4 Proof of main claim 2.3} 
\m
Choose a sufficiently large regular $\l$ and choose a countable $N\prec H_\l$
such that $\rho,P_\d,\le\bar a^\a:\a<\d\r,\tau,M_\d\in N$. Choose a random real
over $N$, $\aa\in(0,1)-\cl(U_{p^*})$. Note that for any rational interval
$(a,b)\sq(0,1)$ we have $\rho((a,b))(\aa)=\chi_{(a,b)}(\aa)$.
\v
Let $u_0=\rho((0,\aa))(\aa)$, $u_1=\rho((\aa,1))(\aa)$. Then $u_0+u_1=1$.
\m
Note that for any number $x$, if $0\leq x<\aa$, then $\rho((x,\aa))(\aa)=u_0$.
[Otherwise, for any rational number $b$ such that $x<b<\aa$, we have
$\rho((0,b))(\aa)>0$, contradicting the choice of $\aa$.] A similar statement
holds for $u_1$. Putting these together we see that $\rho((x,y))(\aa)=1$ for
any numbers $x$ and $y$ such that $0\leq x<\aa<y\leq1$.
\m
Choose an increasing sequence of rational numbers $\le b_{0n}:n<\o\r\in N[\aa]$
converging to $\aa$, and choose a decreasing sequence of rational numbers $\le
b_{1n}:n<\o\r\in N[\aa]$ also converging to $\aa$. In $N[\aa]$ define the
partial order $R$ for adding a Mathias real as follows:
$$R=\{(s,A):s\ \hbox{is a finite subset of}\ \o,\,A\sq\o,\,\max(s)<\min(A)\},$$ 
ordered by $(s,A)\geq(t,B)$
iff $t$ is an initial segment of $s$, $A\sq B$, $s-t\sq B$.
\m
For sets $A\sq\o$, let us identify $A$ with its enumerating function, so that
we may write $A=\{A(i):i<|A|\}$. We need the following special case of the
known fact that an infinite subset of a Mathias real is a Mathias real.  (See
{\bf[M 77:} Theorem 2.0{\bf]}; the
special case which we need here is a fairly routine exercise.)
\m
\n{\bf2.5 Fact} {\it If $X\sq\o$ is $R$-generic over $N[\aa]$, and
$g\in\o^\o\cap N[\aa]$ is increasing, then $Y=\{X(g(n)):n<\o\}$ is also 
$R$-generic over $N[\aa]$.} \endproof
\m
Let $f^*$ be the enumerating function of a set which is  $R$-generic over 
$N[\aa]$.
\m
In $N[\aa][f^*]$, define for increasing functions $f\in\o^\o$,
$$A^k_m(f)=\bigcup_{j<2}\ \bigcup_{k\leq\ell<\o}(b_{j,f(4\ell+m)},
b_{j,f(4\ell+m+1)}).$$
Define $f^*_3(\ell)=f^*(3\ell)$ for $\ell<\o$.
\v
Then $\{A^0_m(f^*_3):m<4\}$ is a partition of $(b_{0,f^*(0)},b_{1,f^*(0)})$. For
some $m<4$ we have 
$$\rho(A^0_m(f^*_3))(\aa)\leq1/4.\leqno(*)$$
\n{\bf2.6 Claim} {\it For any $\bar m<4$ and $k<\o$, we can find an increasing 
function
$g\in N[\aa]\cap\o^\o$ such that $g(i)=i$ for all $i<k$ and
$\rho(A^0_{\bar m}(f^*\circ g))(\aa)\geq3/4$.}
\m
\n{\bf Proof of claim} Let $g(i)=i$ for $i<4k+\bar m+1$ and define 
$g(4\ell+\bar m+1+j)=12\ell+3m+j$ for $\ell\geq k$ and $j<4$. We leave it for
the reader to check, using $(*)$, that $g$ has the desired property. (The reader
might find it helpful, for seeing the role of $g$, to mark off the first few
elements of its range on a line.) \endproof
\m
Let us provisionally let $\bar a^\d=(\le b_{0,f^*(\ell)}:\ell<\o\r,\le
b_{1,f^*(\ell)}:\ell<\o\r,\aa)$. 
\m
\n{\bf2.7 Proof of condition (A) of main claim 2.3} 
\v
Let $J\sq P_\d$ be predense, $J\in M_\d$. We must show that $J$ is predense in
$P_{\d+1}$. Let $p\in P_{\d+1}$, $p\not\in P_\d$. By the definition of
$P_{\d+1}$, there are $q\in P_\d$ and rational numbers $c_0,c_1$ and
$\ell(0)\in\o$ such that 
$$0<b_{0,f^*(4\ell(0))-1}<c_0<b_{0,f^*(4\ell(0))}
<\aa<b_{1,f^*(4\ell(0))}<c_1<b_{1,f^*(4\ell(0))-1}<1,$$
$\cl(U_q)\cap[c_0,c_1]=\e$, $U_p=U_q\cup A^{\ell(0)}_0(f^*)\cup
A^{\ell(0)}_2(f^*)$,
$f_p=f_q\cup0_{A^{\ell(0)}_0(f^*)}\cup1_{A^{\ell(0)}_2(f^*)}$. 
(For $i=0,1$, $i_A$ denotes the function with domain $A$ and constant value
$i$.)
\m
The proof of the following fact is exactly as in \shelah.
\v
\n{\bf2.8 Fact} {\it If $r\in P_\d$, $J\sq P_\d$ is dense, $(c_o,c_1)\sq(0,1)$
and $(c_0,c_1)\cap U_r=\e$, then 
$$\mu\bigl((c_0,c_1)\cap\bigcap\{\cl(U_{r_1}):r_1\in J,r_1\geq r\}\bigr)=0.
\endproof$$} 
Let $J_1=\{r\in P_\d:\exists q_1\in J\ q_1\leq r\}$. For every
$k>f^*(4\ell(0))$ let 
\v
\settabs\+\ $T_k=\{t\in P_\d:$&\cr
\+\ $T_k=\{t\in P_\d:U_t$ is the union of finitely many intervals whose
endpoints are from\cr
\+&$\{b_{j,\ell}:j<2,f^*(4\ell(0))\leq\ell\leq k\}$ and 
$\mu(U_q\cup U_t)<1/2\}$.\cr
\vskip .1cm
So $T_k$ is finite and for each $t\in T_k$, $q\leq q\cup t\in P_\d$ and 
$\aa\not\in\cl(U_t)$. In $N$, define for each $k>f^*(4\ell(0))$ and $t\in T_k$,
\vskip .1cm
\centerline{$J_t=(b_{0,k},b_{1,k})\cap\bigcap\{\cl(U_{r_1}):r_1\in J_1,\,
r_1\geq q\cup t\}$.}
\vskip .1cm
By fact 2.8, $J_t$ has measure zero, and hence $\aa\not\in J_t$.
Thus there is an $r_t\in J_1$, such
that $r_t\geq q\cup t$ and $\aa\not\in\cl(U_{r_t})$.
Define $g\colon\bigcup\{T_k:k>f^*(4\ell(0))\}\to\omega$ and 
$G\colon\omega\to\omega$ such that $[b_{0,g(t)},b_{1,g(t)}]
\cap\cl(U_{r_t})=\emptyset$, $b_{1,g(t)}-b_{0,g(t)}<(1/2)-\mu(U_{r_t})$,
$G(k)=\max\{g(t):t\in T_k\}$.
\v
Since $f^*$ is $R$-generic over $N[\aa]$, for all but finitely many  
$\ell<\o$ we have
$$f^*(4\ell+2)\geq G(f^*(4\ell+1)).$$ 
Choose such an $\ell\geq\ell(0)$.
Let $k=f^*(4\ell+1)$, $t=(U_t,f_t)$, where 
$$U_t=U_p\cap\bigl(
[b_{0,f^*(4\ell(0))},b_{0,k}]\cup[b_{1,k},b_{1,f^*(4\ell(0))}]\bigr),$$
$f_t=f_p|U_t$.
Then $t\in T_k$ and we have $r_t\in J_1$, $r_t\geq q\cup t$. Also, 
$[b_{0,G(k)},b_{1,G(k)}]\cap\cl(U_{r_t})=\emptyset$ and hence
$[b_{0,f^*(4\ell+2)},b_{1,f^*(4\ell+2)}]\cap\cl(U_{r_t})=\emptyset$.
Thus $p$ and $r_t$
are compatible, and this proves part (A) of main claim 2.3.
\m
\n{\bf2.9 Proof of condition (B) of main claim 2.3} 
\m
Let 
$$p^*_1=(U_p\cup A^k_0(f^*)\cup A^k_2(f^*),f_{p^*}\cup0_{A^k_0(f^*)}\cup
1_{A^k_2(f^*)})$$
where $k$ is large enough so that $p_1^*\in P_{\d+1}$. So
$p_1^*\in N[\aa][f^*]$ and $(p^*_1,r^*)\geq(p^*,r^*)$. In $N[\aa][f^*]$, choose
$(p',r')\geq(p_1^*,r^*)$ deciding whether $\tau(\aa)\geq1/2$ or
$\tau(\aa)\leq1/2$, say the first. We will get $(p',r')$ so that condition 
(B1) of main claim 2.3 is
satisfied. The other case is handled similarly. For some $(t,B)\in R\cap
N[\aa]$ we have $f^*(n)=t(n)$ for all $n<|t|$, $f^*(n)\in B$ for all
$n\geq|t|$, and
\m
\centerline{$N[\aa]\models(t,B)\forces_R$``$(p',r')\forces_{P_{\d+1}\times Q}
\,\tau(\aa)\geq1/2$''.} 
\m
\n By claim 2.6 and fact 2.5 above, we can replace $f^*$ by another $R$-generic
real, maintaining $f^*(n)=t(n)$ for $n<|t|$ and $f^*(n)\in B$ for $n\geq|t|$, 
so that $\rho(A^0_2(f^*))(\aa)\geq3/4$. (B1) is
now satisfied.  This completes the proof of main claim 2.3 and of theorem 1.2.
\endproof
\m
\n{\bf Acknowledgment} We thank the referee for several suggestions which
helped improve the readability of the paper.
\vfill\eject
\noindent{\bf References}:
\vskip .2cm
{\bf[BJ 91]} Burke M.R., Just W., {\it Liftings for Haar measure on
$\{0,1\}^\k$}, Israel J. Math, ?? (1991) ???--???.

{\bf[J 89]} Just W., {\it A modification of Shelah's oracle-cc with 
applications}, to appear in Trans. Amer. Math. Soc.

{\bf[F 89]} Fremlin D.H., {\it Measure algebras}, in Handbook of Boolean 
Algebra, ed. J.D. Monk, North-Holland, 1989, 877--980.

{\bf[IT 69]} Ionescu Tulcea A., Ionescu Tulcea C., {\it Topics in the theory of
liftings}, Springer, 1969.

{\bf[M 77]} Mathias A.R.D., {\it Happy families}, Ann. Math. Logic, 12 (1977)
59--111.

{\bf[S 82]} Shelah S., {\it Proper Forcing}, Lecture Notes in Math 940,
Springer-Verlag, 1982.

{\bf[S 83]} Shelah S., {\it Lifting problem of the measure algebra}, Israel J.
Math., 45 (1983) 90--96.
\vskip 1.5cm
\noindent Authors' addresses:
\vskip .3cm
\settabs\+2) &\hfill\cr
\+1) &Department of Mathematics and Computer Science\cr
\+&University of Prince Edward Island\cr
\+&Charlottetown, P.E.I.  C1A 4P3\cr
\+&Canada\cr
\vskip .3cm
\+2) &Institute of Mathematics\cr
\+&The Hebrew University of Jerusalem\cr
\+&Jerusalem, Israel\cr
\end